# BALLISTIC RANDOM WALKS IN RANDOM ENVIRONMENT AT LOW DISORDER

By Christophe Sabot

*CNRS, Université Paris 6*

We consider random walks in a random environment of the type $p_0 + \gamma \xi_z$, where $p_0$ denotes the transition probabilities of a stationary random walk on $\mathbb{Z}^d$, to nearest neighbors, and $\xi_z$ is an i.i.d. random perturbation. We give an explicit expansion, for small $\gamma$, of the asymptotic speed of the random walk under the annealed law, up to order 2. As an application, we construct, in dimension $d \geq 2$, a walk which goes faster than the stationary walk under the mean environment.

Multidimensional random walks in random environment (RWRE) have received a considerable attention these last few years. In particular, several important qualitative results have been obtained, as a law of large numbers, a central limit theorem under certain conditions (by Sznitman, Zerner and Sznitman in the case of independent identically distributed (i.i.d.) environment, [9, 12]; cf. [2, 10] for a review). Unfortunately, in dimension $d \geq 2$, there are, at the present time, very few quantitative results. For example, the Kalikow's condition, under which the law of large numbers is satisfied with a nonnull velocity, is not very explicit (nor the conditions (T), (T′), [10]), and we have very few information about the parameters entering the law of large numbers and the central limit theorem. The aim of this paper is to give an expansion of the asymptotic velocity of the RWRE, which is the parameter entering the law of large numbers, in the case of an environment obtained as a small i.i.d. perturbation of an homogeneous walk.

In this article we consider random walks in random environment in $\mathbb{Z}^d$, for an environment of the type $p_0 + \gamma \xi_z$, where $p_0$ is the transition probabilities of a deterministic, stationary, random walk on $\mathbb{Z}^d$, and $\xi_z$ an i.i.d. random perturbation. We make the assumption that the mean drift, that is, the drift

---

of the mean environment $p_0 + \gamma \mathbb{E}(\xi_z)$, is nonnull for $\gamma$ small enough, $\gamma \neq 0$. In this case, for small $\gamma$'s, the random walk in random environment $X_n$, with transition probabilities $p_0 + \gamma \xi_z$, has a ballistic behavior with speed $v^\gamma \neq 0$, which means that

$$\lim_{n \to \infty} X_n/n = v^\gamma \qquad \text{a.s.}$$

In this article we give an expansion of $v^\gamma$, up to order 2, for small $\gamma$. The second term of the expansion quantifies the first order of the interaction between the randomness of the environment and the global behavior given by the mean environment. From this formula, we can exhibit an acceleration phenomenon, which is specific to dimension $d \geq 2$: we can construct explicitly some environment for which the speed $v^\gamma$ of the RWRE is larger than the speed of the walk under the mean environment: this cannot happen in dimension 1, where there is always slowdown (similar phenomenons are shown in [3]).

The proofs of our results are mainly based on the auxiliary random walk introduced by Kalikow [6], and the result of Sznitman and Zerner on balistic RWRE [12]. In [6], Kalikow gave a formula which expresses the expectation of the Green function of the RWRE (killed at its exit of a bounded domain) as the Green function of an auxiliary random walk. The transition probabilities of this random walk are obtained as weighted expectation of the environment. Under a certain condition, usually named Kalikow's condition, this was used in [6] to prove the transience of the RWRE in dimension $d \geq 2$, and later by Sznitman and Zerner to prove a law of large number, compare [12]. As we show in Section 3, this auxiliary random walk can also give some information about the asymptotic speed of the RWRE. The expansion of $v^\gamma$ is obtained by an expansion of the transition probabilities of the auxiliary random walk. This relies on estimates of perturbed Green functions. These estimates are easy to derive when $p_0$ has a drift, but are much finer when $p_0$ has no drift.

Let us now point out that random perturbations of simple random walks have been considered in several different works (cf. [3, 4, 11]). The type of perturbation we consider here is very specific: the mean drift is at least of order $\gamma$, which is also the order of the perturbation. More precisely, it means that under our assumption (H) (cf. Section 1), the drift of the mean environment at a single point, $\mathbb{E}(p_0 + \gamma \xi_z)$, is at least of order $\gamma$, which is also the order of the perturbation $\gamma \xi_z$. For comparison, when $p_0$ is the transition probability of the simple random walk, then our assumption (H) is stronger than the assumption (0.9) of [11], which implies that the RWRE has a ballistic behavior. When the drift is smaller than the order of the perturbation, different phenomenons may appear as diffusive behavior with nonnull static drift (cf. [3]).

**1. Statement of the results.** In this paper we consider random walks in random environment in $\mathbb{Z}^d$, $d \geq 1$, to nearest neighbors. We denote by $(e_1, \ldots, e_d)$ the canonical basis of the lattices $\mathbb{Z}^d$. We denote by $\mathcal{V}$ the set of elements $e$ in $\mathbb{Z}^d$ with $|e| = 1$, that is, $\mathcal{V} = \{\pm e_1, \ldots, \pm e_d\}$. We will consider environments of the form

$$\omega^\gamma(z, e) = p_0(e) + \gamma \xi(z, e),$$

for $z$ in $\mathbb{Z}^d$ and $e \in \mathcal{V}$, where $p_0$ is a vector of $]0,1[^{\mathcal{V}}$ such that $\sum_{e \in \mathcal{V}} p_0(e) = 1$, and $(\xi(z, \cdot))_{z \in \mathbb{Z}^d}$ are i.i.d. random variables with values in $[-1,1]^{\mathcal{V}}$ and common law $\nu$, and such that $\sum_{e \in \mathcal{V}} \xi(z, e) = 0$, $\nu$-almost surely. We denote by $\mu = \nu^{\otimes \mathbb{Z}^d}$, the law of $(\xi(z, \cdot))_{z \in \mathbb{Z}^d}$. Clearly, there exists $\gamma_0 > 0$ and $\kappa_0 > 0$, such that for all $\gamma < \gamma_0$, we have $\kappa_0 < \omega^\gamma(z, e) < 1$ for all $z$ and $e$, $\mu$-almost surely.

We denote by $P_{z_0}^{\omega^\gamma}$ the law of the Markov chain to nearest neighbors on $\mathbb{Z}^d$, starting from $z_0$, and with transition probability

$$P_{z_0}^{\omega^\gamma}[X_{n+1} = z + e | X_n = z] = \omega^\gamma(z, e) \qquad \forall z, z_0 \in \mathbb{Z}^d, e \in \mathcal{V},$$

and by

$$\mathbb{P}_{z_0}^\gamma(\cdot) = \mathbb{E}_\mu(P_{z_0}^{\omega^\gamma}(\cdot)),$$

the annealed law, where $\mathbb{E}_\mu$ denotes the expectation with respect to $\mu$. The aim of this paper is to give an expansion of the asymptotic speed of the random walk under $\mathbb{P}_{z_0}^\gamma$, up to order 2, in the limit of small $\gamma$'s.

Let us introduce some notation. We set

$$p_1(e) = \mathbb{E}_\mu(\xi(0, e)), \qquad \overline{\xi}(z, e) = \xi(z, e) - p_1(e), \qquad p^\gamma(e) = p_0(e) + \gamma p_1(e)$$

and

$$d_0 = \sum_{e \in \mathcal{V}} p_0(e) e, \qquad d_1 = \sum_{e \in \mathcal{V}} p_1(e) e.$$

We also set

$$C_{e,e'} = \mathrm{Cov}_\nu(\xi(0, e), \xi(0, e')) = \mathbb{E}_\nu(\overline{\xi}(0, e) \overline{\xi}(0, e')).$$

Let us make the following assumption:

(H) $d_0 \neq 0$ or $d_1 \neq 0$.

Then, for $\gamma$ small enough, $d_0 + \gamma d_1 \neq 0$ and the stationary random walk with transition probability $p^\gamma = p_0 + \gamma p_1$ is transient. We denote by $G^{p^\gamma}(z, z')$ the Green function of this walk, and we set

$$J_e^\gamma = G^{p^\gamma}(e, 0) - G^{p^\gamma}(0, 0) \qquad \forall e \in \mathcal{V}.$$

Finally, we set
$$p_{2,\gamma}(e) = \sum_{e' \in \mathcal{V}} C_{e,e'} J^\gamma_{e'}$$
and
$$d_{2,\gamma} = \sum_{e \in \mathcal{V}} p_{2,\gamma}(e) e.$$

THEOREM 1. (i) *For $\gamma$ small enough, $\gamma \neq 0$, $X_n$ is ballistic under $\mathbb{P}^\gamma_\cdot$, that is, there exists $v^\gamma \in \mathbb{R}^d$, $v^\gamma \neq 0$, such that*
$$\lim_{n \to \infty} \frac{X_n}{n} = v^\gamma, \qquad \mathbb{P}^\gamma\text{-}a.s.$$

(ii) *The asymptotic speed $v^\gamma$ has the following expansion for small $\gamma$, $\gamma \neq 0$,*
$$v^\gamma = d_0 + \gamma d_1 + \gamma^2 d_{2,\gamma} + O(\gamma^{3-\varepsilon}),$$
*for all $\varepsilon > 0$.*

REMARK 1. Of course, the interesting part of this formula is the term $d_{2,\gamma}$, which quantifies the interaction between the randomness of the walk, contained in the term $C_{e,e'}$, and the global behavior of the walk, contained in the term $J^\gamma$.

REMARK 2. When $d_0 \neq 0$, we can prove much better estimates, and even get the third order of the expansion. This is the object of Theorem 3 in Section 6.

REMARK 3. This result is valid for $d \geq 1$, in particular, it includes $d = 1$ (where the explicit value of $v^\gamma$ is known). As we shall see, for $d \geq 2$, the second-order term $d_{2,\gamma}$ can be replaced by a term $d_2$, independent of $\gamma$, at the order $O(\gamma^{3-\varepsilon})$. This is not the case for $d = 1$, where $d_{2,\gamma}$ have a discontinuity at $\gamma = 0$ when $d_0 = 0$.

We can prove the previous result by giving an explicit expansion of the term $J^\gamma_e$.

*Dimension* 1. This case is not very interesting since the explicit value of the speed is known [8], but we include it here for completeness. In this case, $J^\gamma_{\pm e_1}$ can be computed explicitly and we have the following: if $(d_0 + \gamma d_1) \cdot e_1 > 0$, then

(1) $$J^\gamma_{-e_1} = 0, \qquad J^\gamma_{e_1} = -\frac{1}{p_0(e_1) + \gamma p^1(e_1)},$$

and the symmetric result holds when $(d_0 + \gamma d_1) \cdot e_1 < 0$. In particular, when $d_0 \neq 0$ and $d_0 \cdot e_1 > 0$, then $J^\gamma_{-e_1} = 0$ and $J^\gamma_{e_1} = -\frac{1}{p_0(e_1)} + O(\gamma)$. It gives that

$$d_{2,\gamma} = -2\frac{1}{p_0(e_1)}\mathbb{E}_\nu(\overline{\xi}(0,e_1)^2)e_1 + O(\gamma).$$

When $d_0 = 0$ [i.e., $p_0(e_1) = p_0(-e_1) = \frac{1}{2}$], and $d_1 \cdot e_1 > 0$, then there is a discontinuity in the second-order term at $\gamma = 0$, and

(2) $$d_{2,\gamma} = -4\operatorname{sgn}(\gamma)\mathbb{E}_\nu(\overline{\xi}(0,e_1)^2)e_1 + O(\gamma).$$

*Dimension $d \geq 2$.* If $d_0 \neq 0$, then we have

(3) $$\begin{aligned}J^\gamma_{\pm e_i} &= \frac{1}{(2\pi)^d}\left(\sqrt{\frac{p_0(\mp e_i)}{p_0(\pm e_i)}} - 1\right) \\ &\quad \times \int_{[0,2\pi]^d} \frac{\cos u_i}{1 - 2\sum_{j=1}^d \sqrt{p_0(e_j)p_0(-e_j)}\cos(u_j)} \prod du_j \\ &\quad + \frac{1}{(2\pi)^d} \int_{[0,2\pi]^d} \frac{(\cos(u_i) - 1)}{1 - 2\sum_{j=1}^d \sqrt{p_0(e_j)p_0(-e_j)}\cos(u_j)} \prod du_j \\ &\quad + O(\gamma).\end{aligned}$$

If $d_0 = 0$, $d_1 \neq 0$, then we have $p_0(e_i) = p_0(-e_i)$ and

(4) $$\begin{aligned}J^\gamma_{\pm e_i} &= \frac{1}{(2\pi)^d}\int_{[0,2\pi]^d}\frac{(\cos(u_i)-1)}{1-2\sum_{j=1}^d p_0(e_j)\cos(u_j)}\prod du_j \\ &\quad + \begin{cases}O(\gamma\log\gamma), & \text{for } d=2, \\ O(\gamma), & \text{for } d \geq 3.\end{cases}\end{aligned}$$

NOTE. Remark that when $d_0 \neq 0$, then $2\sum_{i=1}^d \sqrt{p_0(e_i)p_0(-e_i)} < 1$ and the integrand in both terms of (3) is bounded. When $d_0 = 0$, the integrand in (4) is also bounded due to the presence of the term $(\cos(u_i) - 1)$ in the numerator.

In these two cases, we write $J_e$ for the first term of the expansion of $J^\gamma_e$, which is independent of $\gamma$. Hence, we see that for $d \geq 2$, the expansion of Theorem 1 can be rewritten

(5) $$v^\gamma = d_0 + \gamma d_1 + \gamma^2 d_2 + O(\gamma^{3-\varepsilon}),$$

with $d_2 = \sum_{e \in \mathcal{V}} p_2(e)e$, where $p_2(e) = \sum_{e'} C_{e,e'} J_{e'}$.

REMARK 4. In dimension $d = 2$, the second term of the expansion of $J^\gamma_e$ induces the Green function of a symmetric walk killed at rate $K\gamma^2$, for some $K > 0$, compare Section 4. This Green function diverges like $\log \gamma$, and the estimate we give in Theorem 1 does not allow to include this term in the expansion of $v^\gamma$. We think our estimates in Theorem 1 could be improved in order to allow us to include this term.

*Heuristic interpretation.* Remark that the term $\gamma^2 d_{2,\gamma}$ can also be written

$$\gamma^2 d_{2,\gamma} = \mathbb{E}_\nu\left[\left(\sum_{e \in \mathcal{V}} \gamma \overline{\xi}(0,e)e\right)\left(\gamma \sum_{e \in \mathcal{V}} \overline{\xi}(0,e)(G^{p^\gamma}(e,0) - G^{p^\gamma}(0,0))\right)\right]$$
(6)
$$= \mathbb{E}_\nu\left[\left(\sum_{e \in \mathcal{V}} \gamma \overline{\xi}(0,e)e\right)\left(1 - \frac{G^{p^\gamma}(0,0)}{G^{\tilde{p}^\gamma}(0,0)}\right)\right],$$

where $\tilde{p}^\gamma$ is the one point modification of $p^\gamma$, given by $\tilde{p}^\gamma(z,\cdot) = p^\gamma(z,\cdot)$, $z \neq 0$, $\tilde{p}^\gamma(0,\cdot) = \omega^\gamma(0,\cdot)$. Hence, we see that $\gamma^2 d_{2,\gamma}$ is a weighted expectation of the random part of the drift $\gamma \sum_e \overline{\xi}(0,e)e$. The weight is positive when $G^{\tilde{p}^\gamma}(0,0) > G^{p^\gamma}(0,0)$, that is, when the statistical number of visit of 0 is increased by the randomization of the environment at the point 0. This means that at this order, the environment has a larger weight when the point is more often visited. In fact, this is one of the pieces of information contained in Kalikow's formula (cf. Section 3). The interest of this formula is to quantify this effect.

Let us now explain the structure of the paper. In Section 2 we apply these results to show that the speed $v^\gamma$ can be larger than the speed of the stationary walk under the mean environment. In Section 3 we recall the definition of the auxiliary random walk introduced by Kalikow and the law of large numbers proved by Sznitman and Zerner, and give a simple result relating the effective value of the asymptotic speed with the drift of the Kalikow's random walk. In Section 4 we prove Theorem 1. Section 5 is devoted to the proof of the formulas concerning $J_e^\gamma$. In Section 6 we give the third order of the expansion, when $d_0 \neq 0$.

**2. Speedup in higher dimension.** Considering the formula of $d_{2,\gamma}$ in the case $d = 1$, we see that the second-order drift is in opposite direction to the main drift $d_0 + \gamma d_1$. It is actually true for any balistic RWRE in dimension 1, that the asymptotic speed of the RWRE is smaller than the mean drift given by the random environment. Indeed, in dimension 1, if we consider a RWRE with i.i.d. environment $w_z$, then this walk has a ballistic behavior in the positive direction if and only if $\mathbb{E}(\frac{\omega(-e_1)}{\omega(e_1)}) < 1$, and in this case, the asymptotic speed has the following expression (cf. [8]):

$$\frac{1 - \mathbb{E}(\omega(-e_1))/\omega(e_1)}{1 + \mathbb{E}(\omega(-e_1)/\omega(e_1))} = \frac{\mathbb{E}((\omega(e_1) - \omega(-e_1))/\mathbb{E}(\omega(e_1)))}{\mathbb{E}(1/\omega(e_1))},$$

which is easily seen to be smaller than $\mathbb{E}(\omega(e_1) - \omega(-e_1))$. The intuitive explanation for this slowdown effect is that the sites where the environment plays against the main behavior are overweighted, in the sense that the expected number of visits of these sites is larger.

This phenomenon is no longer valid in higher dimension. We construct here an explicit RWRE, for which the asymptotic speed is larger than the mean drift at one site.

REMARK 5. Similar phenomenons are shown in [3], where the authors exhibit RWRE with positive velocity and negative mean drift in large dimension (cf. Remark 4.3.2 and Section 5.3 of [3]).

REMARK 6. We give here an example in dimension $d = 2$, but we could easily give a similar example in any dimension $d \geq 2$.

Let us consider $d = 2$, and $p_0$ given by

$$p_0(e_1) = p_0(-e_1) = \frac{1+a}{4},$$

$$p_0(e_2) = \frac{1-a}{4}(1+\varepsilon),$$

$$p_0(-e_2) = \frac{1-a}{4}(1-\varepsilon),$$

for some reals $0 < a < 1$, and $0 < \varepsilon < 1$. We see that

$$d_0 = \varepsilon \frac{1-a}{2} e_2.$$

Let us now define $U \in \mathbb{R}^{\mathcal{V}}$ by

$$U(e_1) = U(-e_1) = 1,$$
$$U(e_2) = 0, \qquad U(-e_2) = -2,$$

and the random variable $(\xi(z, \cdot))$ by

$$\xi(z, \cdot) = U(\cdot) \qquad \text{with probability } \tfrac{1}{2},$$
$$\xi(z, \cdot) = -U(\cdot) \qquad \text{with probability } \tfrac{1}{2},$$

independently on each site $z$ in $\mathbb{Z}^2$. It is clear that $\mathbb{E}(\xi(z, e)) = 0$ for all $z, e$ and, hence, that $p_1 = 0$. The covariance matrix is given by

$$\text{Cov}(\xi(\pm e_1), \xi(\pm e_1)) = 1,$$
$$\text{Cov}(\xi(\pm e_1), \xi(e_2)) = \text{Cov}(\xi(\pm e_2), \xi(e_2)) = 0,$$
$$\text{Cov}(\xi(-e_2), \xi(-e_2)) = 4,$$
$$\text{Cov}(\xi(\pm e_1), \xi(-e_2)) = -2,$$

where we simply wrote $\xi(e)$ for $\xi(z,e)$. It is clear, by symmetry, that $d_2$ will be in the direction $\pm e_2$ [we recall that $d_2$ is defined in (5)], and computation gives

$$d_2 = (2(J_{-e_1} + J_{e_1}) - 4J_{-e_2})e_2,$$

where $J_{\pm e_i}$ is the first order in $\gamma$ of $J_e^\gamma$, given by (3). When $\varepsilon$ goes to zero, the first term in formula (3) goes to zero [indeed, the term $(\sqrt{\frac{p_0(-e_2)}{p_0(e_2)}} - 1)$ is of order $\varepsilon$ and the second term is of order $\log \varepsilon$, since it is the Green function at 0 of a stationary Markov chain with killing rate of order $\varepsilon^2$ (cf. the discussion of Section 4)]. It implies that

$$\lim_{\varepsilon \to 0} d_2 = \left( \int_{[0,2\pi]^2} \frac{(\cos(u_1) - 1) - (\cos(u_2) - 1)}{1 - (1/2)((1-a)\cos(u_2) + (1+a)\cos(u_1))} \right) e_2$$
$$= \left( \int_{[0,2\pi]^2} \frac{2(\cos(u_1) - \cos(u_2))}{(1 - \cos(u_2)) + (1 - \cos(u_1)) + a(\cos(u_2) - \cos(u_1))} \right) e_2.$$

It is not difficult to check that the previous integral is positive. Indeed, if we consider $U_1$ and $U_2$, two uniform random variables on $[0, 2\pi]$, and $S = \cos(U_1) + \cos(U_2)$, $A = \cos(U_2) - \cos(U_1)$, then, by symmetry, we have $\mathbb{E}(A|S) = 0$, and since the previous integral is equal to

$$\mathbb{E}\left(\frac{-2A}{2 - S + aA}\right) = \mathbb{E}\left(\mathbb{E}\left(\frac{-2A}{2 - S + aA}\bigg|S\right)\right),$$

we classically get that it is positive, when $a$ is positive. This implies that for $\varepsilon$ small enough, the term $d_2$ is in the direction $+e_2$, hence, in the same direction as $d_0$. It implies that for $\gamma$ small enough, $v^\gamma \cdot e_2 > d_0 \cdot e_2$.

REMARK 7. The intuitive explanation for this phenomena is that, due to the nonsymmetry of the horizontal and vertical direction, it is easier for the walk, under the mean environment $p_0$, to come back to 0 from the point $\pm e_1$ than from the point $\pm e_2$. Then we choose a random environment which correlates the acceleration of the walk, that is, a drift in the direction $e_2$ larger than the mean drift, with a larger probability to go on the horizontal direction. This overweights, in Kalikow's formula, the environment which have a larger drift in the direction $e_2$ than the mean drift. Indeed, in formula (6), for $\varepsilon$ small enough, the drift $\gamma \sum_e \overline{\xi}(0,e)e$ has a positive weight when $\xi = U$, and a negative weight when $\xi = -U$ (indeed, the one-point modification of $p^\gamma$ in $\tilde{p}^\gamma$ increases the statistical number of visit of 0, when $\xi = U$, and decreases it when $\xi = -U$). But, $(\gamma \sum_e \overline{\xi}(0,e)e) \cdot e_2$ is positive when $\xi = U$ and negative when $\xi = -U$.

**3. Kalikow's auxiliary random walk.** We present a generalization of the random walk introduced by Kalikow in [6].

Let us first introduce some notation. Let $\kappa_0 > 0$ be a positive real. We denote by $\Omega_{\kappa_0}$ the set of environments with uniform ellipticity constant $\kappa_0$, that is,

$$\Omega_{\kappa_0} = \bigg\{ (w(x,e))_{x \in \mathbb{Z}^d, e \in \mathcal{V}} \in (\,]0,1[^{2d})^{\mathbb{Z}^d},$$

$$\text{such that } \sum_{e \in \mathcal{V}} \omega(x,e) = 1 \text{ and } \omega(x,e) \geq \kappa_0 \; \forall x, e \bigg\}.$$

We suppose that $\mu$ is a probability measure on $\Omega_{\kappa_0}$. (We do not assume, for the moment, that $\mu$ is the law of an i.i.d. environment.) Let $U$ be a connected subset of $\mathbb{Z}^d$, and $\delta$ a real $0 < \delta \leq 1$. We denote by $\partial U$ the boundary set of $U$, that is, $\partial U = \{z \in \mathbb{Z}^d \setminus U, \; \exists x \in U, |z - x| = 1\}$. If $\omega \in \Omega$ is an environment, we denote, as usual, by $P_z^\omega$ the law of the random walk in the environment $\omega$, starting from $z$. We set, for $z \in U$, $z' \in U \cup \partial U$,

$$G_{U,\delta}^\omega(z,z') = E_z^\omega\bigg(\sum_{k=0}^{T_U} \delta^k \mathbb{1}_{\{X_k = z'\}}\bigg),$$

where $T_U = \inf\{k, X_k \in \mathbb{Z}^d \setminus U\}$. For $z' \in U$, $G_{U,\delta}^\omega(z,z')$ is the value of the usual Green function of the process killed at constant rate $\delta$, and stopped after its first hitting time of $\mathbb{Z}^d \setminus U$, and for $z' \in \partial U$,

$$G_U^\omega(z,z') = E_z^\omega(\mathbb{1}_{\{X_{T_U} = z'\}} \delta^{T_U})$$

is the probability to exit $U$ at the point $z'$, before having been killed.

Let us now fix a point $z_0$ in $U$. We suppose that either $U$ is bounded, or $\delta < 1$. For all $z$ in $U$, we set

$$\hat{w}_{U,\delta,z_0}(z,e) = \frac{\mathbb{E}_\mu(G_{U,\delta}^\omega(z_0,z)\omega(z,e))}{\mathbb{E}_\mu(G_{U,\delta}^\omega(z_0,z))}.$$

Obviously, $(\hat{\omega}_{U,\delta,z_0}(z,\cdot))$ defines the transition probabilities of a Markov chain on $U$. We denote by $G_U^{\hat{\omega}_{U,\delta,z_0}}$ the Green function of this Markov chain, killed at constant rate $\delta$, and stopped after the first exit time of $U$,

$$G_{U,\delta}^{\hat{\omega}_{U,\delta,z_0}}(z,z') = E_z^{\hat{\omega}_{U,\delta,z_0}}\bigg(\sum_{k=0}^{T_U} \delta^k \mathbb{1}_{\{X_k = z'\}}\bigg).$$

We simply write $G_U^\omega$ and $\hat{\omega}_{U,z_0}$, when $\delta = 1$ and $U$ is a bounded subset of $\mathbb{Z}^d$, and $G_\delta$, $\hat{\omega}_{\delta,z_0}$, when $U = \mathbb{Z}^d$ and $\delta < 1$.

In its generalized version, the result of Kalikow says the following:

PROPOSITION 1. *If either $U$ is bounded or $\delta < 1$, then for all $z$ in $U \cap \partial U$,*

$$\mathbb{E}_\mu(G_{U,\delta}^\omega(z_0, z)) = G_{U,\delta}^{\hat{\omega}_{U,\delta,z_0}}(z_0, z).$$

REMARK 8. The original result of Kalikow was given for $U$ bounded and $\delta = 1$.

REMARK 9. In the sequel, Kalikow's formula will refer to the formula of Proposition 1.

PROOF. The proof is essentially the same as the proof of Kalikow, but we give it for convenience since the hypothesis are not exactly the same. Let us first remark that $0 < G_{U,\delta}^\omega(x,y) < c_{x,y}$ for a constant independent of $\omega$ in $\Omega_{\kappa_0}$ (indeed, this comes from the uniform ellipticity condition). This implies that $\hat{\omega}_{U,\delta,z_0}$ is well defined. Remark now that for all $z$ in $U \cup \partial U$,

$$(7) \qquad G_{U,\delta}^\omega(z_0, z) = \delta_{z_0, z} + \sum_{\substack{e \in \mathcal{V}, \text{ s.t.} \\ z-e \in U}} G_{U,\delta}^\omega(z_0, z-e) \delta \omega(z-e, e),$$

which gives

$$\mathbb{E}_\mu(G_{U,\delta}^\omega(z_0, z)) = \delta_{z_0, z} + \sum_{\substack{e \in \mathcal{V}, \text{ s.t.} \\ z-e \in U}} \mathbb{E}_\mu(G_{U,\delta}^\omega(z_0, z-e)) \delta \hat{\omega}_{U,\delta,z_0}(z-e, e).$$

Let us set

$$G_{U,\delta}^{(n)}(z_0, z) = E_{z_0}^{\hat{\omega}_{U,\delta,z_0}}\left(\sum_{k=0}^{n \wedge T_U} \delta^k \mathbb{1}_{\{X_k = z'\}}\right).$$

We see that

$$G_{U,\delta}^{(n+1)}(z_0, z) = \delta_{z_0, z} + \sum_{\substack{e \in \mathcal{V}, \text{ s.t.} \\ z-e \in U}} G_{U,\delta}^{(n)}(z_0, z-e) \delta \hat{\omega}_{U,\delta,z_0}(z-e, e).$$

It is clear that

$$G_{U,\delta}^{\hat{\omega}_{U,\delta,z_0}}(z_0, z) = \lim_{n \to \infty} G_{U,\delta}^{(n)}(z_0, z),$$

and that $G_{U,\delta}^{\hat{\omega}_{U,\delta,z_0}}$ satisfies the same equation as $G_{U,\delta}^\omega$ in (7). By induction, we have

$$G_{U,\delta}^{(n)}(z_0, z) \leq \mathbb{E}_\mu(G_{U,\delta}^\omega(z_0, z))$$

for all $n$ and, thus,

$$(8) \qquad G_{U,\delta}^{\hat{\omega}_{U,\delta,z_0}}(z_0, z) \leq \mathbb{E}_\mu(G_{U,\delta}^\omega(z_0, z))$$

for all $z$ in $U \cup \partial U$. For $\delta < 1$, it is clear that for all environment $\omega$ in $\Omega_{\kappa_0}$,

$$\frac{1}{1-\delta} = \sum_{k=0}^{\infty} \delta^k$$

$$= E_{z_0}^{\omega}\left(\sum_{k=0}^{T_U-1} \delta^k + \delta^{T_u}\frac{1}{1-\delta}\right)$$

$$= E_{z_0}^{\omega}\left(\sum_{z \in U}\sum_{k=0}^{T_u-1} \delta^k \mathbb{1}_{\{X_k=z\}} + \frac{1}{1-\delta}\sum_{z \in \partial U} \delta^{T_U}\mathbb{1}_{\{X_{T_U}=z\}}\right)$$

$$= \sum_{z \in U} G_{U,\delta}^{\omega}(z_0, z) + \frac{1}{1-\delta}\sum_{z \in \partial U} G_{U,\delta}^{\omega}(z_0, z).$$

This implies both

$$\frac{1}{1-\delta} = \sum_{z \in U} \mathbb{E}_\mu(G_{U,\delta}^{\omega}(z_0, z)) + \frac{1}{1-\delta}\sum_{z \in \partial U} \mathbb{E}_\mu(G_{U,\delta}^{\omega}(z_0, z))$$

and

$$\frac{1}{1-\delta} = \sum_{z \in U} G_{U,\delta}^{\hat{\omega}_{U,\delta,z_0}}(z_0, z) + \frac{1}{1-\delta}\sum_{z \in \partial U} G_{U,\delta}^{\hat{\omega}_{U,\delta,z_0}}(z_0, z).$$

This necessarily implies equality in (8) for all $z$ in $U \cup \partial U$.

When $\delta = 1$ and $U$ is bounded, then $T_U < \infty$ almost surely under $P_{z_0}^\omega$, for any elliptic environment $\omega$. This implies

$$1 = \sum_{z \in \partial U} G_U^{\hat{\omega}_{U,z_0}}(z_0, z) = \sum_{z \in \partial U} \mathbb{E}_\mu(G_U^\omega(z_0, z)),$$

which gives equality in (8) for all $z$ in $\partial U$, and then by induction for all $z$ in $U$.

□

3.1. *Application to the asymptotic speed.* We first recall a result of Sznitman and Zerner (cf. [12]). Let us suppose now that $\mu = \nu^{\otimes \mathbb{Z}^d}$ is the law of a uniformly elliptic, i.i.d. random environment $(w(x,\cdot))$ in $\Omega_{\kappa_0}$.

THEOREM 2 [12]. *If there exists a vector $l$ in $\mathbb{R}^d$, a constant $\varepsilon > 0$, such that for all bounded connected subset $U \subset \mathbb{Z}^d$, $U \neq \varnothing$, and all $z_0 \in U$,*

$$\left(\sum_{e \in \mathcal{V}} \hat{\omega}_{U,z_0}(z,e)e\right) \cdot l \geq \varepsilon \qquad \forall z \in U,$$

*then there exists a vector $v \in \mathbb{R}^d$, such that*

$$\lim_{n \to \infty} \frac{X_n}{n} = v, \qquad P^\omega\text{-}a.s.,$$

*for $\mu$-almost all environment $\omega$. Moreover, $v \cdot l > 0$, hence, $X_n$ is ballistic in the direction $l$.*

REMARK 10. The velocity can be expressed in terms of the expectation of some renewal time (cf. [12]) or in terms of Lyapounov exponents [13], but we will not need these expressions.

We suppose now that our RWRE satisfies the condition of the previous theorem. We can easily get some information on the asymptotic speed from the walk of Kalikow. We consider now $U = \mathbb{Z}^d$ and $\delta < 1$. The transition probabilities, $\hat{w}_{\delta,z_0}(z,e)$, of the Kalikow's walk depend only on the difference $z - z_0$. We denote by $\hat{d}_\delta(z) = \sum_{e \in \mathcal{V}} \hat{\omega}_{\delta,0}(z,e) e$ the drift associated with the Kalikow's walk. We denote by $\mathcal{A}_\delta$ the convex hull of the set

$$\bigcup_{z \in \mathbb{Z}^d} d_\delta(z),$$

and by $\mathcal{A}$, the set of accumulation points of $\mathcal{A}_\delta$, when $\delta$ goes to 1.

PROPOSITION 2. *The asymptotic speed $v$ is in $\mathcal{A}$.*

PROOF. We denote by $\mathbb{E}_0$ the expectation with respect to the annealed law

$$\mathbb{E}_0(\cdot) = \mathbb{E}_\mu(E_0^\omega(\cdot)).$$

We consider an independent geometric random variable $\tau_\delta$ with parameter $\delta$. We have

$$\begin{aligned}
\mathbb{E}_0(X_{\tau_\delta}) &= \mathbb{E}_\mu\left(E_0^\omega\left(\sum_{k=0}^{\tau_\delta - 1} X_{k+1} - X_k\right)\right) \\
&= \sum_{z \in \mathbb{Z}^d} \sum_{k=0}^\infty \mathbb{E}_\mu(E_0^\omega(\mathbb{1}_{\{\tau_\delta > k, x_k = z\}}(X_{k+1} - X_k))) \\
&= \sum_{z \in \mathbb{Z}^d} \sum_{k=0}^\infty \mathbb{E}_\mu\left(\mathbb{1}_{\{\tau_\delta > k, x_k = z\}}\left(\sum_{e \in \mathcal{V}} \omega(z,e) e\right)\right) \\
&= \sum_{z \in \mathbb{Z}^d} \mathbb{E}_\mu(G_\delta^\omega(0,z))\left(\sum_{e \in \mathcal{V}} \hat{\omega}_{\delta,0}(z,e) e\right) \\
&= \sum_{z \in \mathbb{Z}^d} G_\delta^{\hat{\omega}_{\delta,0}}(0,z) \hat{d}_\delta(z).
\end{aligned}$$

But $\mathbb{E}(\tau_\delta) = \sum_{z \in \mathbb{Z}^d} \hat{G}_\delta^{\hat{\omega}_{\delta,0}}(0,z)$. Thus,

$$\frac{\mathbb{E}_0(X_{\tau_\delta})}{\mathbb{E}(\tau_\delta)} = \frac{\sum_{z \in \mathbb{Z}^d} G_\delta^{\hat{\omega}_{\delta,0}}(0,z) \hat{d}_\delta(z)}{\sum_{z \in \mathbb{Z}^d} G_\delta^{\hat{\omega}_{\delta,0}}(0,z)} \in \mathcal{A}_\delta,$$

and since $\lim_{\delta \to 1} \mathbb{E}_0(X_{\tau_\delta})/\mathbb{E}(\tau_\delta) = v$, we know that $v$ is in $\mathcal{A}$. □

**4. Proof of Theorem 1.** We will use the following simple estimates several times in the article.

LEMMA 1. *Let $\omega$ and $\omega'$ be two environments in $\Omega_{\kappa_0}$ for some $\kappa_0 > 0$. We suppose that $\omega'$ is a perturbation of $\omega$, at some point $z$ in $\mathbb{Z}^d$, that is, that we have*

$$\omega'(z',e) = \omega(z',e) \qquad \text{for } z' \neq z, e \in \mathcal{V},$$
$$\omega'(z,e) = \omega(z,e) + \Delta\omega(e) \qquad \text{for } e \in \mathcal{V},$$

*for some $(\Delta\omega(e)) \in \,]-1,1[^\mathcal{V}$. Let $U \subset \mathbb{Z}^d$, $0 < \delta \leq 1$, be such that either $U$ is bounded or $\delta < 1$, and such that $z$ is in $U$. Then, we have the following estimates: for all $y$ in $U$, $y'$ in $U \cup \partial U$,*

$$|G_{U,\delta}^{\omega'}(y,y') - G_{U,\delta}^{\omega}(y,y')| \leq \frac{2d \sup_{e \in \mathcal{V}} |\Delta\omega(e)|}{\kappa_0^2} G_{U,\delta}^{\omega'}(y,y')$$

*and*

$$\left| G_{U,\delta}^{\omega'}(y,y') - G_{U,\delta}^{\omega}(y,y') - G_{U,\delta}^{\omega}(y,z) \sum_{e \in \mathcal{V}} \Delta\omega(e)(\delta G_{U,\delta}^{\omega}(z+e,y') - G_{U,\delta}^{\omega}(z,y')) \right|$$
$$\leq \frac{(2d \sup_{e \in \mathcal{V}} |\Delta\omega(e)|)^2}{\kappa_0^3} G_{U,\delta}^{\omega'}(y,y').$$

PROOF. To simplify notation, in this proof we simply write $G^\omega$ for $G_{U,\delta}^\omega$. We will use several times in this paper the following classical expansion of Green functions: let $P$ and $P'$ be the transition operators of two random walks on $\mathbb{Z}^d$ (with eventually some killing, so that $0 \leq \sum_y P_{x,y} \leq 1$, the left-hand side inequality being eventually strict), to nearest neighbors, and $G_\delta^P = (I - \delta P)^{-1} = \sum \delta^k P^k$, and $G_\delta^{P'} = (I - \delta P')^{-1}$ the associated Green functions (for $0 < \delta < 1$). Then for all $n \geq 0$, we have

$$(9) \quad G_\delta^{P'} = G_\delta^P + \sum_{k=1}^{n} \delta^k (G_\delta^P (P' - P))^k G_\delta^P + \delta^{n+1} (G_\delta^P (P' - P))^{n+1} G_\delta^{P'}.$$

In particular, for $n=0$, we have

(10) $$G_\delta^{P'} = G_\delta^P + \delta G_\delta^P (P' - P) G_\delta^{P'}.$$

We apply the previous formula for the transition operators associated with transition probabilities $\omega$ and $\omega'$, and we get

(11) $$\begin{aligned} G^{\omega'}(y, y') &- G^\omega(y, y') \\ &= \delta G^\omega(y, z) \sum_{e \in \mathcal{V}} \Delta\omega(e) G^{\omega'}(z+e, y') \\ &= G^\omega(y, z) \sum_{e \in \mathcal{V}} \Delta\omega(e)(\delta G^{\omega'}(z+e, y') - G^{\omega'}(z, y')). \end{aligned}$$

In the last formula, we used that $\sum_{e \in \mathcal{V}} \Delta\omega(e) = 0$. If we set $T_z = \inf\{n \geq 0, X_n = z\}$, we get

(12) $$\begin{aligned} \delta G^{\omega'}&(z+e, y') - G^{\omega'}(z, y') \\ &\geq (\delta \mathbb{E}_{z+e}^{\omega'}(\mathbb{1}_{\{T_z < T_U\}} \delta^{T_z}) - 1) G^{\omega'}(z, y') \\ &= (\delta \mathbb{E}_{z+e}^{\omega}(\mathbb{1}_{\{T_z < T_U\}} \delta^{T_z}) - 1) G^{\omega'}(z, y') \\ &\geq \frac{1}{\kappa_0} G^{\omega'}(z, y') \sum_{e' \in \mathcal{V}} \omega(z, e')(\delta \mathbb{E}_{z+e'}^{\omega}(\mathbb{1}_{\{T_z < T_U\}} \delta^{T_z}) - 1) \\ &= -\frac{1}{\kappa_0} \frac{G^{\omega'}(z, y')}{G^\omega(z, z)}. \end{aligned}$$

Remark that with the same argument, we also have

$$\delta G^{\omega'}(z+e, y') - G^{\omega'}(z, y') \geq -\frac{1}{\kappa_0} \frac{G^{\omega'}(z, y')}{G^{\omega'}(z, z)},$$

which is equal to $-\frac{1}{\kappa_0}$, if $y' = z$. If $y' = z$, then $\delta G^{\omega'}(z+e, z) - G^{\omega'}(z, z) \leq 0$. In particular, this gives, when $y' = z$,

(13) $$|\delta G^{\omega'}(z+e, z) - G^{\omega'}(z, z)| \leq \frac{1}{\kappa_0}.$$

If $z \neq y'$, then

$$\delta \sum_{e' \in \mathcal{V}} \omega'(z, e') G^{\omega'}(z+e', y') = G^{\omega'}(z, y'),$$

which gives

$$\begin{aligned} (\delta G^{\omega'}&(z+e, y') - G^{\omega'}(z, y')) \\ &= -\sum_{e' \in \mathcal{V}, e' \neq e} \frac{\omega'(z, e')}{\omega'(z, e)} (\delta G^{\omega'}(z+e', y') - G^{\omega'}(z, y')) \\ &\leq \frac{1}{\kappa_0^2} \frac{G^{\omega'}(z, y')}{G^\omega(z, z)}, \end{aligned}$$

where in the last inequality we used estimate (12). Thus, we get

$$(14) \qquad |\delta G^{\omega'}(z+e,y') - G^{\omega'}(z,y')| \leq \frac{1}{\kappa_0^2} \frac{G^{\omega'}(z,y')}{G^{\omega}(z,z)}.$$

Applied to (11), it gives

$$|G^{\omega'}(y,y') - G^{\omega}(y,y')| \leq \frac{2d \sup_{e \in \mathcal{V}} |\Delta \omega(e)|}{\kappa_0^2} \mathbb{E}_y^{\omega}(\mathbb{1}_{\{T_z < T_U\}} \delta^{T_z}) G^{\omega'}(z,y')$$

$$\leq \frac{2d \sup_{e \in \mathcal{V}} |\Delta \omega(e)|}{\kappa_0^2} G^{\omega'}(y,y').$$

The second estimate is similar. We expand $G^{\omega'}$ at order 2 [i.e., we use (9) with $n=1$] which gives

$$G^{\omega'}(y,y') - G^{\omega}(y,y') - \sum_{e \in \mathcal{V}} G^{\omega}(y,z) \Delta \omega(e) (\delta G^{\omega}(z+e,y') - G^{\omega}(z,y'))$$

$$= \sum_{e \in \mathcal{V}} \sum_{e' \in \mathcal{V}} G^{\omega}(y,z) \Delta \omega(e) (\delta G^{\omega}(z+e,z) - G^{\omega}(z,z))$$

$$\times \Delta \omega(e') (\delta G^{\omega'}(z+e',y') - G^{\omega'}(z,y')).$$

But, $|\delta G^{\omega}(z+e,z) - G_\delta^{\omega}(z,z)| \leq \frac{1}{\kappa_0}$, compare (13), and using (14), we get the second estimate. □

In order to prove Theorem 1, we give an expansion of Kalikow's transition probabilities. This is based on two successive applications of Proposition 1. We come back to the notation of Section 1. We have an i.i.d. random environment of the form

$$\omega^{\gamma}(x,e) = p_0(e) + \gamma \xi(x,e) = p_0(e) + \gamma(p_1(e) + \overline{\xi}(x,e)),$$

where $\xi(x,e)$ is distributed according to the law $\mu = \nu^{\otimes \mathbb{Z}^d}$, and

$$p_1(e) = \mathbb{E}_\mu(\xi(x,e)), \qquad \overline{\xi}(x,e) = \xi(x,e) - p_1(e) \qquad \forall x \in \mathbb{Z}^d, e \in \mathcal{V}.$$

For any $y$ in $\mathbb{Z}^d$, we denote by $\omega^{\gamma,y}$ the environment

$$\omega^{\gamma,y}(z,e) = \begin{cases} \omega^{\gamma}(z,e), & \text{if } z \neq y, \\ p_0(e) + \gamma p_1(e), & \text{if } z = y. \end{cases}$$

For $0 \leq \delta \leq 1$ and $U \subset \mathbb{Z}^d$, with either $\delta < 1$ or $U$ bounded, $z_0 \in U$, we denote by $\hat{\omega}_{U,\delta,z_0}^{\gamma}$ the transition probabilities of the auxiliary random walk defined in Section 3, associated with the environment $\omega^{\gamma}$ under $\mu$.

LEMMA 2. *We have the following expansion, for small $\gamma$'s,*

$$\hat{\omega}_{U,\delta,z_0}^{\gamma}(y,e) = p_0(e) + \gamma p_1(e) + \gamma^2 \sum_{e' \in \mathcal{V}} C_{e,e'} \widetilde{J}_{e'}^{U,\delta,z_0,\gamma}(y) + O(\gamma^3),$$

*where*

$$\widetilde{J}_{e'}^{U,\delta,z_0,\gamma}(y) = \frac{\mathbb{E}_\mu(G_{U,\delta}^{\omega^{\gamma,y}}(z_0,y)(\delta G_{U,\delta}^{\omega^{\gamma,y}}(y+e',y) - G_{U,\delta}^{\omega^{\gamma,y}}(y,y)))}{\mathbb{E}_\mu(G_{U,\delta}^{\omega^{\gamma,y}}(z_0,y))},$$

*and where* $|O(\gamma^3)| \leq 2\frac{(2d)^2}{\kappa_0^4}\gamma^3$.

PROOF. We simply write $G^\omega$ for $G_{U,\delta}^\omega$ in this proof. Let us first remark that we have

$$\hat{\omega}_{U,\delta,z_0}^\gamma(y,e) = p^\gamma(e) + \gamma \frac{\mathbb{E}_\mu(G^{\omega^\gamma}(z_0,y)\overline{\xi}(y,e))}{\mathbb{E}_\mu(G^{\omega^\gamma}(z_0,y))},$$

where $p^\gamma = p_0 + \gamma p_1$. Applying Lemma 1 to $\omega^\gamma$ and $\omega^{\gamma,y}$, we get

$$\frac{\mathbb{E}_\mu(G^{\omega^\gamma}(z_0,y)\overline{\xi}(y,e))}{\mathbb{E}_\mu(G^{\omega^\gamma}(z_0,y))}$$

$$= \frac{\mathbb{E}_\mu(G^{\omega^{\gamma,y}}(z_0,y)\overline{\xi}(y,e))}{\mathbb{E}_\mu(G^{\omega^\gamma}(z_0,y))}$$

$$+ \gamma \sum_{e'\in\mathcal{V}} \mathbb{E}_\mu(G^{\omega^{\gamma,y}}(z_0,y)\overline{\xi}(y,e')$$

$$\times (\delta G^{\omega^{\gamma,y}}(y+e',y) - G^{\omega^{\gamma,y}}(y,y))\overline{\xi}(y,e))$$

$$\times [\mathbb{E}_\mu(G^{\omega^\gamma}(z_0,y))]^{-1} + O_1(\gamma^2)$$

$$= \gamma \sum_{e'\in\mathcal{V}} \mathbb{E}_\mu(\overline{\xi}(y,e)\overline{\xi}(y,e'))$$

$$\times \frac{\mathbb{E}_\mu(G^{\omega^{\gamma,y}}(z_0,y)(\delta G^{\omega^{\gamma,y}}(y+e',y) - G^{\omega^{\gamma,y}}(y,y)))}{\mathbb{E}_\mu(G^{\omega^\gamma}(z_0,y))} + O_1(\gamma^2),$$

where $|O_1(\gamma^2)| \leq \frac{(2d)^2}{\kappa_0^3}\gamma^2$. [In the last formula, we used the independence of $G^{\omega^{\gamma,y}}$ and $\overline{\xi}(y,e)$, and the fact that $\mathbb{E}_\mu(\overline{\xi}(y,e)) = 0$.] Considering now that by Lemma 1 we have

$$\left|1 - \frac{\mathbb{E}_\mu(G_\delta^{\omega^{\gamma,y}}(z_0,y))}{\mathbb{E}_\mu(G_\delta^{\omega^\gamma}(z_0,y))}\right| \leq \frac{(2d)^2}{\kappa_0^3}\gamma^2,$$

we get

$$\frac{\mathbb{E}_\mu(G^{\omega^{\gamma,y}}(z_0,y)(\delta G^{\omega^{\gamma,y}}(y+e',y) - G^{\omega^{\gamma,y}}(y,y)))}{\mathbb{E}_\mu(G^{\omega^\gamma}(z_0,y))}$$

$$= \frac{\mathbb{E}_\mu(G^{\omega^{\gamma,y}}(z_0,y)(\delta G^{\omega^{\gamma,y}}(y+e',y) - G^{\omega^{\gamma,y}}(y,y)))}{\mathbb{E}_\mu(G^{\omega^{\gamma,y}}(z_0,y))} + O_2(\gamma^2),$$

where $|O_2(\gamma^2)| \leq \frac{(2d)^2}{\kappa_0^4}\gamma^2$. □

Let us remark that the previous lemma implies that, under the hypothesis (H), for $\gamma$ small enough, $\gamma \neq 0$, there exists a positive constant $c_\gamma$, such that for all bounded connected subset $U$ and $z_0$ in $U$,

$$\left(\sum_{e \in \mathcal{V}} \hat{\omega}_{U,z_0}(z,e)e\right) \cdot (d_0 + \gamma d_1) > c_\gamma \qquad \forall z \in U.$$

[Indeed, $|J_e^{U,\delta,z_0,\gamma}(y)| \leq \frac{1}{\kappa_0}$, using (13).] Hence, we are in the condition of application of Theorem 2, for $l = d_0 + \gamma d_1$. To obtain information on the speed $v^\gamma$, we have to estimate the transition probabilities $\omega_{\delta,0}^\gamma$, when $\delta$ goes to 1. This is the object of the next lemma. We simply write $\widetilde{J}_e^{\delta,\gamma}(y)$ for $\widetilde{J}_e^{U,\delta,z_0,\gamma}(y)$ when $U = \mathbb{Z}^d$ and $z_0 = 0$.

LEMMA 3. (i) If $d_0 \neq 0$, then there exists a constant $C > 0$, such that for $\gamma$ sufficiently small,

$$\limsup_{\delta \to 1} |\widetilde{J}_e^{\delta,\gamma}(y) - J_e^\gamma| \leq C\gamma^2,$$

for all $e \in \mathcal{V}$, $y \in \mathbb{Z}^d$.

(ii) If $d_0 = 0$ and $d_1 \neq 0$, then, for all $0 < \varepsilon < 1$, there exists a constant $C_\varepsilon > 0$, such that for $\gamma$ sufficiently small,

$$\limsup_{\delta \to 1} |\widetilde{J}_e^{\delta,\gamma}(y) - J_e^\gamma| \leq C_\varepsilon \gamma^{1-\varepsilon},$$

for all $e \in \mathcal{V}$, $y \in \mathbb{Z}^d$.

Let us first point out that this lemma concludes the proof of Theorem 1 using Lemma 2 and Proposition 2.

PROOF OF LEMMA 3. Let us first describe the structure of the proof. In the first step, we apply Proposition 1 to a certain modified measure $\tilde{\mu}^y$ to write $\widetilde{J}_e^{\delta,\gamma}(y)$ as the Green function $\delta G_\delta^{p^\gamma + \gamma^2 \Delta \omega}(y+e,y) - G_\delta^{p^\gamma + \gamma^2 \Delta \omega}(y,y)$, for a certain deterministic environment $p^\gamma + \gamma^2 \Delta \omega$, which is a second-order perturbation of $p^\gamma$. In the second step, we expand the Green function $G_\delta^{p^\gamma + \gamma^2 \Delta \omega}$, and rewrite the Green function of $G^{p^\gamma}$ as the Green function of a symmetric random walk plus a killing. The last step is to use estimates on quantities like

$$\sum_{z \in \mathbb{Z}^d} |p_{n+1}(0,z) - p_n(e,z)|,$$

which are adapted from Lawler's book [7].

*Step* 1. For $y$, $y'$, $y''$ in $\mathbb{Z}^d$, we write

$$I^{\delta,\gamma,y}(y',y'') = \frac{\mathbb{E}_\mu(G_\delta^{\omega^{\gamma,y}}(0,y)G_\delta^{\omega^{\gamma,y}}(y',y''))}{\mathbb{E}_\mu(G_\delta^{\omega^{\gamma,y}}(0,y))},$$

so that we have $\widetilde{J}_e^{\delta,\gamma}(y) = \delta I^{\delta,\gamma,y}(y+e,y) - I^{\delta,\gamma,y}(y,y)$. We denote by $\tilde{\mu}^y$ the probability measure on $\Omega_{\kappa_0}$ given by

$$\tilde{\mu}^y = \frac{G_\delta^{\omega^{\gamma,y}}(0,y)}{\mathbb{E}_\mu(G_\delta^{\omega^{\gamma,y}}(0,y))}\mu.$$

It is clear that

$$I^{\delta,\gamma,y}(y',y'') = \mathbb{E}_{\tilde{\mu}^y}(G_\delta^{\omega^{\gamma,y}}(y',y'')).$$

We apply Proposition 1 for the environment $\omega^{\gamma,y}$ under the measure $\tilde{\mu}^y$, for the initial point $z_0 = y'$. This means that the Kalikow's random walk has transition probabilities

$$\tilde{\omega}(z,e) = \frac{\mathbb{E}_{\tilde{\mu}^y}(G_\delta^{\omega^{\gamma,y}}(y',z)\omega^{\gamma,y}(z,e))}{\mathbb{E}_{\tilde{\mu}^y}(G_\delta^{\omega^{\gamma,y}}(y',z))},$$

and that by Proposition 1 we get

$$I^{\delta,\gamma,y}(y',y'') = G_\delta^{\tilde{\omega}}(y',y'').$$

We have

$$\tilde{\omega}(z,e) = p^\gamma(e) + \gamma\frac{\mathbb{E}_{\tilde{\mu}^y}(G_\delta^{\omega^{\gamma,y}}(y',z)\overline{\xi}(z,e))}{\mathbb{E}_{\tilde{\mu}^y}(G_\delta^{\omega^{\gamma,y}}(y',z))},$$

for $z \neq y$, and $\tilde{\omega}(y,e) = p^\gamma(e)$. We want to prove that

$$\tilde{\omega}(z,e) = p^\gamma(e) + \gamma^2 \Delta\omega(z,e),$$

for a perturbative term $\Delta\omega(z,e)$, uniformly bounded in $y,y',z,e,\delta,\gamma$. So we write

$$\Delta\omega(z,e) = \frac{1}{\gamma^2}(\tilde{\omega}(z,e) - p^\gamma(e))$$

for $\gamma \neq 0$. As in Lemma 2, we define the environment $(\omega^{\gamma,y,z}(z',e))_{z',e}$ in $\Omega_{\kappa_0}$, by

$$\omega^{\gamma,y,z}(z',e) = \begin{cases} p^\gamma(e), & \text{if } z' = z \text{ or } z' = y, \\ \omega^\gamma(z',e), & \text{if } z' \neq z,\ z' \neq y. \end{cases}$$

Using Lemma 1, we get

$$\frac{\mathbb{E}_{\tilde{\mu}^y}(G_\delta^{\omega^{\gamma,y}}(y,z)\overline{\xi}(z,e))}{\mathbb{E}_{\tilde{\mu}^y}(G_\delta^{\omega^{\gamma,y}}(y,z))} = \frac{\mathbb{E}_{\tilde{\mu}^y}(G_\delta^{\omega^{\gamma,y,z}}(y,z)\overline{\xi}(z,e))}{\mathbb{E}_{\tilde{\mu}^y}(G_\delta^{\omega^{\gamma,y,z}}(y,z))} + O(\gamma)$$

BALLISTIC RWRE AT LOW DISORDER 19$$= \frac{\mathbb{E}_\mu(G_\delta^{\omega^{\gamma,y}}(0,y)G_\delta^{\omega^{\gamma,y,z}}(y,z)\overline{\xi}(z,e))}{\mathbb{E}_\mu(G_\delta^{\omega^{\gamma,y}}(0,y)G_\delta^{\omega^{\gamma,y,z}}(y,z))} + O(\gamma)$$

$$= \frac{\mathbb{E}_\mu(G_\delta^{\omega^{\gamma,y,z}}(0,y)G_\delta^{\omega^{\gamma,y,z}}(y,z)\overline{\xi}(z,e))}{\mathbb{E}_\mu(G_\delta^{\omega^{\gamma,y,z}}(0,y)G_\delta^{\omega^{\gamma,y,z}}(y,z))} + O(\gamma)$$

$$= O(\gamma),$$

where, as usual, the remainder terms $O(\gamma)$ satisfy $|O(\gamma)| \le C|\gamma|$, where $C > 0$ is a constant depending only on $\kappa_0$, $d$. [In the last equality, we used, as usual, the independence of $G_\delta^{\omega^{\gamma,y,z}}$ with $\overline{\xi}(z,e)$, and the fact that $\mathbb{E}_\mu(\overline{\xi}(z,e)) = 0$.] This implies that $\Delta\omega(z,e)$ is bounded by a constant depending only on $\kappa_0$, $d$.

*Step* 2. We transform now the Green function $G_\delta^{p^\gamma}$ into the Green function of a symmetric walk plus a killing. Let $\phi^\gamma : \mathbb{Z}^d \to \mathbb{R}$ be defined by

$$\phi^\gamma(z) = \prod_{i=1}^d \left(\sqrt{\frac{p^\gamma(e_i)}{p^\gamma(-e_i)}}\right)^{z_i}.$$

Let $M_\phi$ be the operator of multiplication by $\phi$, given, for $f : \mathbb{Z}^d \to \mathbb{R}$, by

$$M_{\phi^\gamma}(f)(z) = \phi^\gamma(z)f(z).$$

If $P^{p^\gamma}$ is the transition operator of the walk with stationary transition probabilities $(p^\gamma(e))_{e \in \mathcal{V}}$, then we have

$$M_{\phi^\gamma} P^{p^\gamma} M_{\phi^\gamma}^{-1} = k^\gamma P^{s^\gamma},$$

where

$$k^\gamma = 2\sum_{i=1}^d (p^\gamma(-e_i)p^\gamma(e_i))^{1/2},$$

and $P^{s^\gamma}$ is the transition operator of the symmetric, stationary random walk, with transition probabilities

$$s^\gamma(e_i) = s^\gamma(-e_i) = \frac{(p^\gamma(-e_i)p^\gamma(e_i))^{1/2}}{k^\gamma}.$$

As we shall see later, $k^\gamma < 1$, and we have

$$G_\delta^{p^\gamma} = (I - \delta P^{p^\gamma})^{-1} = M_{\phi^\gamma}^{-1}(I - \delta k^\gamma P^{s^\gamma})^{-1} M_{\phi^\gamma} = M_{\phi^\gamma}^{-1} G_{\delta k^\gamma}^{s^\gamma} M_{\phi^\gamma}.$$

Let us come back to $k^\gamma$. We trivially have

$$1 - k^\gamma = \sum_{i=1}^d (\sqrt{p^\gamma(e_i)} - \sqrt{p^\gamma(-e_i)})^2,$$

which implies that $k^\gamma < 1$ under the hypothesis (H) for small $\gamma$. If $d_0 \neq 0$, then

$$1 - k^\gamma = \sum_{i=1}^{d} (\sqrt{p_0(e_i)} - \sqrt{p_0(-e_i)})^2 + O(\gamma).$$

If $d_0 = 0$ and $d_1 \neq 0$, then we easily get

$$1 - k^\gamma = \frac{\gamma^2}{4} \sum_{i=1}^{d} \frac{(\sqrt{p_1(e_i)} - \sqrt{p_1(-e_i)})^2}{p_0(e_i)} + O(\gamma^3).$$

It means that in this case $1 - k^\gamma = K\gamma^2 + O(\gamma^3)$ for a positive constant $K > 0$.

We can easily get, similarly, that

$$G_\delta^{p^\gamma + \gamma^2 \Delta\omega} = M_{\phi^\gamma}^{-1} G_{\delta \tilde{k}^\gamma}^{\tilde{s}^\gamma} M_{\phi^\gamma},$$

where $\tilde{s}^\gamma = \tilde{s}^\gamma(z, e)$ is, a priori, a nonsymmetric, nonstationary environment of the form

$$\tilde{s}^\gamma(z, e) = s^\gamma(e) + \gamma^2 \Delta s(z, e),$$

where $\Delta s$ is uniformly bounded independently of the variables $y, y', z, e, \gamma, \delta$. When $d_0 \neq 0$, the term $\tilde{k}^\gamma = \tilde{k}^\gamma(z)$ is of the form

$$\tilde{k}^\gamma(z) = k^\gamma + \gamma^2 \Delta k(z),$$

where $\Delta k$ is uniformly bounded. When $d_0 = 0$, $d_1 \neq 0$, then the term $\tilde{k}^\gamma(z)$ is of the form

$$\tilde{k}^\gamma(z) = k^\gamma + \gamma^3 \Delta k(z),$$

where $\Delta k$ is uniformly bounded [this comes from the fact that $\sum_e \Delta\omega(e) = 0$, and that $\sqrt{\frac{p^\gamma(\pm e_i)}{p^\gamma(\mp e_i)}} - 1 = O(\gamma)$, from which the term of order 2 is null].

*Step* 3. We consider now the following expansion at order $n$ [which is a consequence of the classical expansion of Green functions, cf. (9)]:

$$G_\delta^{p^\gamma + \gamma^2 \Delta\omega}(z, z') - G_\delta^{p^\gamma}(z, z')$$
$$= \sum_{k=1}^{n} (\delta\gamma^2)^k S_k(z, z') + (\delta\gamma^2)^{n+1} R_n(z, z'),$$

where

$$S_n(z, z') = \sum_{z_1,\ldots,z_n} \sum_{e_1,\ldots,e_n} G_\delta^{p^\gamma}(z, z_1) \Delta\omega(z_1, e_1) G_\delta^{p^\gamma}(z_1 + e_1, z_2) \cdots$$
$$\times \Delta\omega(z_n, e_n) G_\delta^{p^\gamma}(z_n + e_n, z')$$

and
$$R_n = \sum_{z'' \in \mathbb{Z}^d} S_n(z, z'') \sum_{e'' \in \mathcal{V}} \Delta\omega(z'', e'') G_\delta^{p^\gamma + \gamma^2 \Delta\omega}(z'' + e'', z').$$

Considering the transformation of step 2, we get

$$S_n(z, z')$$
$$= \phi^\gamma(z' - z) \sum_{\substack{z_1, \ldots, z_n \\ e_1, \ldots, e_n}} G_{\delta k^\gamma}^{s^\gamma}(z, z_1) \Delta\omega(z_1, e_1) \phi^\gamma(-e_1) G_{\delta k^\gamma}^{s^\gamma}(z_1 + e_1, z_2) \cdots$$
$$\times \Delta\omega(z_n, e_n) \phi^\gamma(-e_n) G_{\delta k^\gamma}^{s^\gamma}(z_n + e_n, z')$$

and

$$R_n(z, z')$$
$$= \phi^\gamma(z' - z) \sum_{\substack{z_1, \ldots, z_n, z'' \\ e_1, \ldots, e_n, e''}} G_{\delta k^\gamma}^{s^\gamma}(z, z_1) \Delta\omega(z_1, e_1) \phi^\gamma(-e_1) G_{\delta k^\gamma}^{s^\gamma}(z_1 + e_1, z_2) \cdots$$
$$\times \Delta\omega(z_n, e_n) \phi^\gamma(-e_n) G_{\delta k^\gamma}^{s^\gamma}(z_n + e_n, z')$$
$$\times \Delta\omega(z'', e'') \phi^\gamma(-e'') G_{\delta \tilde{k}^\gamma}^{\tilde{s}^\gamma}(z'' + e'', z').$$

If $d_0 \neq 0$, then

$$|1 - \delta k^\gamma| \leq 2 \sum_{i=1}^d (\sqrt{p_0(e_i)} - \sqrt{p_0(-e_i)})^2 + O(\gamma).$$

Thus, we get [since $\phi^\gamma(e) \leq \frac{1}{\kappa_0}$, for $\gamma < \gamma_0$]

$$S_n(z, z') \leq \phi^\gamma(z' - z) \left( \frac{1}{\kappa_0} 2d \left( \sup_{z,e} |\Delta\omega(z, e)| \right) \sum_{z \in \mathbb{Z}^d} G_{\delta k^\gamma}^{s^\gamma}(0, z) \right)^{n+1}$$
$$= \phi^\gamma(z' - z) \left( \frac{1}{\kappa_0} 2d \left( \sup_{z,e} |\Delta\omega(z, e)| \right) \frac{1}{1 - \delta k^\gamma} \right)^{n+1}$$
$$\leq \phi^\gamma(z' - z) C^{n+1},$$

for some positive constant $C$, depending only on $\kappa_0$, $d$, $p_0$. We can get a similar estimate for the remaining term $R_n(z, z')$ considering that $1 - \tilde{k}^\gamma(z) \sim 1 - \delta k^\gamma$. This implies that for $\gamma$ small enough, the series $\sum_{k=0}^\infty (\delta k^\gamma)^k S_k(z', z)$ is convergent and that

$$G_\delta^{p^\gamma + \gamma^2 \Delta\omega}(z, z') - G_\delta^{p^\gamma}(z, z') = \sum_{k=1}^\infty (\delta\gamma^2)^k S_k(z, z')$$
$$= \phi^\gamma(z' - z) O(\gamma^2).$$

Considering the discussion of step 1, this concludes Lemma 3(i).

If $d_0 = 0$ and $d_1 \neq 0$, then we rewrite

$$S_n(z,z')$$
$$= \phi^\gamma(z'-z) \sum_{\substack{z_1,\ldots,z_n \\ e_1,\ldots,e_n}} G^{s^\gamma}_{\delta k^\gamma}(z,z_1) \Delta\omega(z_1,e_1)$$
$$\times (\phi^\gamma(-e_1) G^{s^\gamma}_{\delta k^\gamma}(z_1+e_1,z_2) - G^{s^\gamma}_{\delta k^\gamma}(z_1,z_2)) \cdots$$
$$\times \Delta\omega(z_n,e_n)(\phi^\gamma(-e_n) G^{s^\gamma}_{\delta k^\gamma}(z_n+e_n,z') - G^{s^\gamma}_{\delta k^\gamma}(z_n,z'))$$
$$\leq \phi^\gamma(z'-z)(2d)^n G^{s^\gamma}_{\delta k^\gamma}(0,0) \left( \sup_{e\in\mathcal{V}} \left| \sum_{z\in\mathbb{Z}^d} \phi^\gamma(-e) G^{s^\gamma}_{\delta k^\gamma}(e,z) - G^{s^\gamma}_{\delta k^\gamma}(0,z) \right| \right)^n$$
$$\leq \phi^\gamma(z'-z)(2d)^n G^{s^\gamma}_{\delta k^\gamma}(0,0)$$
$$\times \sup_{e\in\mathcal{V}} \left( \left| \frac{\phi^\gamma(-e)-1}{1-\delta k^\gamma} \right| + \left| \sum_{z\in\mathbb{Z}^d} G^{s^\gamma}_{\delta k^\gamma}(e,z) - G^{s^\gamma}_{\delta k^\gamma}(0,z) \right| \right)^n.$$

We write $p_n(y,y')$ for the $n$th step transition probability of the random walk with transition probability $s^\gamma$. We consider the term

$$\left| \sum_{z\in\mathbb{Z}^d} G^{s^\gamma}_{\delta k^\gamma}(e,z) - G^{s^\gamma}_{\delta k^\gamma}(0,z) \right|$$
$$= \left| \sum_{z\in\mathbb{Z}^d} \sum_{n\in\mathbb{N}} (\delta k^\gamma)^n (p_n(e,z) - p_n(0,z)) \right|$$
$$\leq 1 + \left| \sum_{z\in\mathbb{Z}^d} \sum_{n\in\mathbb{N}} (\delta k^\gamma)^n p_n(e,z) - (\delta k^\gamma)^{n+1} p_{n+1}(0,z) \right|$$
$$\leq 1 + (1-\delta k^\gamma) \left| \sum_{z\in\mathbb{Z}^d} \sum_{n\in\mathbb{N}} (\delta k^\gamma)^n p_n(e,z) \right|$$
$$\quad + \sum_{n\in\mathbb{N}} (\delta k^\gamma)^n \left| \sum_{z\in\mathbb{Z}^d} p_n(e,z) - p_{n+1}(0,z) \right|$$
$$\leq 2 + \sum_{n\in\mathbb{N}} (\delta k^\gamma)^n \left| \sum_{z\in\mathbb{Z}^d} p_n(e,z) - p_{n+1}(0,z) \right|.$$

We now use the following lemma, which specifies and generalizes Corollary 1.2.3 of [7] (we will prove this lemma later on).

LEMMA 4. *Let* $(s(e))_{e\in\mathcal{V}} \in ]0,1[^{\mathcal{V}}$ *be such that*

$$s(e_i) = s(-e_i) \geq \kappa_0, \qquad 2\sum_{i=1}^{d} s(e_i) = 1.$$

*Then, for all $\varepsilon > 0$, there exists a positive constant $C_\varepsilon$, depending only on $\kappa_0$, $d$, such that*

$$\sum_{z\in\mathbb{Z}^d} |p_{n+1}(0,z) - p_n(e,z)| \leq C_\varepsilon n^{-(1/2)+\varepsilon} \qquad \forall e \in \mathcal{V},$$

*where $p_n(z,z')$ is the $n$th step transition probability of the stationary, symmetric, random walk on $\mathbb{Z}^d$, with transition probability $(s(e))_{e\in\mathcal{V}}$.*

NOTE. The difference in the indices $n$ and $n+1$ comes from the fact that $p_n(z,z')$ is null if $n$ and $\sum_j z'_j - z_j$ do not have the same parity.

It means that for all $0 < \varepsilon < 1$, there is a positive constant $C_\varepsilon > 0$, such that

$$\left|\sum_{z\in\mathbb{Z}^d} G^{s^\gamma}_{\delta k^\gamma}(e,z) - G^{s^\gamma}_{\delta k^\gamma}(0,z)\right| \leq 2 + C_\varepsilon \sum_{n=0}^{\infty} (\delta k^\gamma)^n n^{-(1/2)+\varepsilon}.$$

Considering that for $x > 0$,

$$(1-\delta k^\gamma)^{-x} = \sum_{n=0}^{\infty} \frac{x(x+1)\cdots(x+n-1)}{n!}(\delta k^\gamma)^n,$$

and that $\frac{x(x+1)\cdots(x+n-1)}{n!} \asymp n^{x-1}$ for large $n$ [which means that there is a constant $\widetilde{K} > 0$, for which $\widetilde{K}^{-1}n^{x-1} \leq \frac{x(x+1)\cdots(x+n-1)}{n!} \leq \widetilde{K}n^{x-1}$], we see that for all $x$ such that $(1-x) < \frac{1}{2} - \varepsilon$, we can find a constant $C > 0$ such that

$$\sum_{n=0}^{\infty} (\delta k^\gamma)^n n^{-(1/2)+\varepsilon} \leq C(1-\delta k^\gamma)^{-x}.$$

This means that for all $\varepsilon > 0$, we can find a new constant $C_\varepsilon > 0$, such that

$$\left|\sum_{z\in\mathbb{Z}^d} G^{s^\gamma}_{\delta k^\gamma}(e,z) - G^{s^\gamma}_{\delta k^\gamma}(0,z)\right| \leq C_\varepsilon(1-\delta k^\gamma)^{-(1/2)-\varepsilon}$$

$$\leq C_\varepsilon(1-\delta k^\gamma)^{-(1/2)-\varepsilon}.$$

Considering that $1 - k^\gamma \sim K\gamma^2$ for small $\gamma$'s, we see that, for all positive $\varepsilon$, we can find a new constant $C_\varepsilon > 0$, such that

$$\left|\sum_{z\in\mathbb{Z}^d} G^{s^\gamma}_{\delta k^\gamma}(e,z) - G^{s^\gamma}_{\delta k^\gamma}(0,z)\right| \leq C_\varepsilon \gamma^{-1-\varepsilon}.$$

Coming back to $S_n(z,z')$ and considering that $\phi^\gamma(-e) - 1 = O(\gamma)$, we see that

$$|S_n(z,z')| \leq \phi^\gamma(z'-z)(O(\gamma^{-1-\varepsilon}))^n |G^{s^\gamma}_{\delta k^\gamma}(0,0)|$$
$$\leq \phi^\gamma(z'-z)(O(\gamma^{-1-\varepsilon}))^n |G^{s^\gamma}_{k^\gamma}(0,0)|.$$

It clearly implies that $\sum_{k\geq 1}(\delta\gamma^2)^k S_k(z,z')$ is absolutely convergent for $\gamma$ sufficiently small, and that

$$\sum_n (\delta\gamma^2)^n S_n(z,z') = \phi^\gamma(z'-z) O(G^{s^\gamma}_{k^\gamma}(0,0)\gamma^{1-\varepsilon}),$$

for all positive $\varepsilon$. It is not difficult, using the same arguments and the fact that $1 - \tilde{k}^\gamma(z) = K\gamma^2 + O(\gamma^3)$ for the same constant $K > 0$, to prove that the remaining term $R_n$ goes to 0, when $n$ goes to infinity, which means that

$$G^{p^\gamma}_\delta(z,z') - G^{p^\gamma + \gamma^2 \Delta\omega}_\delta(z,z') = \phi^\gamma(z-z') \times O(G^{s^\gamma}_{k^\gamma}(0,0)\gamma^{1-\varepsilon})$$

and the estimate $O(G^{s^\gamma}_{k^\gamma}(0,0)\gamma^{1-\varepsilon})$ is uniform in $\delta$, $z$, $z'$. Coming back to step 1, we see that it means that

$$\limsup_{\delta \to 0} |\widetilde{J}^{\delta,\gamma}_e(y) - J^\gamma_e| = O(G^{s^\gamma}_{k^\gamma}(0,0)\gamma^{1-\varepsilon}).$$

Since, in dim 2, $G^{s^\gamma}_{k^\gamma}(0,0)$ diverges like $\ln(1 - k^\gamma) \sim \ln\gamma$, and is bounded in dimension $d \geq 3$, we see that it means that

$$\limsup_{\delta \to 1} |\widetilde{J}^{\delta,\gamma}_e(y) - J^\gamma_e| = O(\gamma^{1-\varepsilon})$$

for all positive $\varepsilon$. This concludes the proof for $d \geq 2$. For $d = 1$, the expansion of $v^\gamma$ can be checked directly from the explicit formula for $v^\gamma$. (It could also be deduced from the same method.) □

PROOF OF LEMMA 4. Let $s_i = s(e_i) = s(-e_i)$. We suppose that $n+1$ and $\sum_{j=1}^d z_j$ have the same parity, since, otherwise, $p_{n+1}(0,z)$ and $p_n(e,z)$ are null. By the Fourier transform we have

$$p_{n+1}(0,z) - p_n(\pm e_i, z)$$
$$= \frac{1}{(2\pi)^d} \int_{[0,2\pi]^d} \Theta^{n+1} \prod_{j=1}^d \cos(z_j u_j)\, du_j$$
$$- \frac{1}{(2\pi)^d} \int_{[0,2\pi]^d} \Theta^n \cos((z_i \mp 1)u_i) \prod_{j\neq i} \cos(z_j u_j) \prod_{j=1}^d du_j,$$

where

$$\Theta = \sum_{j=1}^d s_j \cos(u_j).$$

We can find $\rho < 1$, $r > 0$, such that $|\Theta| \leq \rho$ if $(u_1, \ldots, u_d) \notin\, ]-r, r[^d \cup\, ]\pi - r, \pi + r[^d$. We take constants $C > 0$ and $0 < \eta < 1$ such that

$$|1 - \cos u| \leq Cu, \qquad |\sin u| \leq C|u|, \qquad |\cos u| \leq \left|1 - \frac{(\eta u)^2}{2}\right| \qquad \forall u \in\, ]-r, r[.$$

Hence, we have, using parity of the integrands,

$$|p_{n+1}(0, z) - p_n(\pm e_i, z)|$$

$$\leq 2\rho^n + \frac{2}{(2\pi)^d} \int_{]-r,r[^d} \left(\sum_{j=1}^d s_j(\cos(u_j) - 1)\right) \Theta^n \prod_{j=1}^d \cos(z_j u_j)\, du_j$$

$$+ \frac{2}{(2\pi)^d} \int_{]-r,r[^d} (1 - \cos u_i) \Theta^n \prod_{j=1}^d \cos(z_j u_j)\, du_j$$

$$+ \frac{2}{(2\pi)^d} \int_{]-r,r[^d} \Theta^n \sin u_i \sin z_i u_i \prod_{j \neq i} \cos(z_j u_j)\, du_j$$

$$\leq 2\rho^n + C'|z_i| \int_{]-r,r[^d} \left(\sum_{j=1}^d u_j^2\right) \left(1 - \frac{1}{2} \sum_{j=1}^d s_j(\eta u_j)^2\right)^n \prod_{j=1}^d du_j,$$

with $C' = \frac{6C}{(2\pi)^d}$. Considering $\alpha_j = \sqrt{n} u_j$, the last expression is equal to

$$2\rho^n + C' \frac{1}{n^{(d+2)/2}} |z_i| \int_{|\alpha_j| \leq r\sqrt{n}} \left(\sum_{j=1}^d \alpha_j^2\right) \left(1 - \frac{1}{2} \sum_{j=1}^d s_j(\eta \alpha_j)^2/n\right)^n \prod_{j=1}^d d\alpha_j.$$

Considering the inequality $(1 - \frac{1}{2} u^2/n)^n \leq \exp(-u^2/2)$ for $u \leq \sqrt{2n}$, we get that the last expression is smaller than

$$2\rho^n + C' \frac{1}{n^{(d+2)/2}} |z_i| \int_{\mathbb{R}^d} \left(\sum_{j=1}^d \alpha_j^2\right) \exp\left(1 - \frac{1}{2} \sum_{j=1}^d s_j(\eta \alpha_j)^2\right) \prod_{j=1}^d d\alpha_j$$

$$= |z_i| O\left(\frac{1}{n^{(d+2)/2}}\right).$$

[Let us remark that this $O(\frac{1}{n^{(d+2)/2}})$ can be uniformly estimated on the set of transition probabilities $s$ satisfying the uniform ellipticity condition with constant $\kappa_0$.]

Let us take $\theta > \frac{1}{2}$. For $e \in \mathcal{V}$, we consider the sum

$$\sum_{z \in \mathbb{Z}^d} |p_{n+1}(0, z) - p_n(e, z)|$$

$$= \sum_{|z| > n^\theta} |p_{n+1}(0, z) - p_n(e, z)| + \sum_{|z| \leq n^\theta} |p_{n+1}(0, z) - p_n(e, z)|.$$

The first term is bounded by $4\exp-\frac{1}{2d}n^{2\theta-1}$ using Hoeffding's inequality (cf., e.g., [5]). For the second term we use the previous estimate

$$\sum_{|z|\leq n^\theta} |p_{n+1}(0,z) - p_n(e,z)| \leq O\left(\frac{1}{n^{(d+2)/2}}\right) \sum_{|z|\leq n^\theta} |z|$$

$$= O(n^{\theta(d+1)-(d+2)/2})$$

$$= O(n^{\theta-1+d(\theta-(1/2))}).$$

Since we can take any $\theta > \frac{1}{2}$, we can get any order $O(n^{-(1/2)+\varepsilon})$ for $\varepsilon > 0$. □

## 5. Development of $J_e^\gamma$.

*In dimension $d=1$.* When $(d_0 + \gamma d_1) \cdot e_1 > 0$, the random walk with stationary probability $p^\gamma$ is transcient in the positive direction. This implies that

$$P_{-e_1}^{p^\gamma}(T_0 < \infty) = 1$$

and, thus, that

$$G^{p^\gamma}(-e_1, 0) - G^{p^\gamma}(0,0) = 0.$$

Hence,

$$1 = -p^\gamma(e_1)(G^{p^\gamma}(e_1, 0) - G^{p^\gamma}(0,0)),$$

which gives formula (1).

*Dimension $d \geq 2$.* We use the transformation described in step 2 of the proof of Lemma 3 to get

$$G^{p^\gamma}(z,z') = \phi^\gamma(z'-z)G_{k^\gamma}^{s^\gamma}(z,z').$$

It gives

$$J_e^\gamma = (\phi^\gamma(-e) - 1)G_{k^\gamma}^{s^\gamma}(e,0) + (G_{k^\gamma}^{s^\gamma}(e,0) - G_{k^\gamma}^{s^\gamma}(0,0)).$$

Using the Fourier transform and since $k^\gamma s^\gamma(\pm e_i) = \sqrt{p^\gamma(e_i)p^\gamma(-e_i)}$, we get

$$J_{\pm e_i}^\gamma = \frac{1}{(2\pi)^d}\left(\sqrt{\frac{p_\gamma(\mp e_i)}{p_\gamma(\pm e_i)}} - 1\right)$$

$$\times \int_{[0,2\pi]^d} \frac{\cos u_i}{1 - 2\sum_{j=1}^d \sqrt{p^\gamma(e_j)p^\gamma(-e_j)}\cos(u_j)} \prod du_j$$

$$+ \frac{1}{(2\pi)^d}\int_{[0,2\pi]^d} \frac{\cos(u_i) - 1}{1 - 2\sum_{j=1}^d \sqrt{p^\gamma(e_j)p^\gamma(-e_j)}\cos(u_j)} \prod du_j.$$

When $d_0 \neq 0$, it is clear that the last formula gives (3), at first order in $\gamma$ [since $2\sum_{j=1}^{d}\sqrt{p_0(e_j)p_0(-e_j)} < 1$, which implies that the denominator is uniformly bounded away from 0]. When $d = 2$ and $d_0 = 0$, $d_1 \neq 0$, then $(\phi^\gamma(-e) - 1)G_{k\gamma}^{s\gamma}(e,0)$ is of order $\gamma \log \gamma$. When $d_0 = 0$, $d_1 \neq 0$, $d \geq 3$, $G_{k\gamma}^{s\gamma}(e,0)$ is uniformly bounded, and the first term is of order $O(\gamma)$. In any case, the second term is uniformly bounded and gives (4) at first order in $\gamma$.

**6. The third order when $d_0 \neq 0$.** When $d_0 \neq 0$, we can improve Theorem 1.

THEOREM 3. *If $d_0 \neq 0$, then*
$$v^\gamma = d_0 + \gamma d_1 + \gamma^2 d_{2,\gamma} + \gamma^3 d_3 + O(\gamma^4),$$
*where $d_3 = \sum_{e \in \mathcal{V}} p_3(e)e$, and*
$$p_3(e) = \sum_{e',e'' \in \mathcal{V}} \mathbb{E}_\mu(\overline{\xi}(e)\overline{\xi}(e')\overline{\xi}(e''))J_{e'}J_{e''}.$$

NOTE. $J_e$ is the first order of the expansion of $J_e^\gamma$, compare (3).

PROOF. We just sketch the proof, since it is simple and very similar to the proof of the second order in the case $d_0 \neq 0$. We can improve Lemma 2 as follows (we only consider the case $\delta < 1$ and $U = \mathbb{Z}^d$):
$$\hat{\omega}_{\delta,0}^\gamma(y,e) = p_0(e) + \gamma p_1(e) + \gamma^2 \sum_{e' \in \mathcal{V}} \mathbb{E}_\mu(\overline{\xi}(e)\overline{\xi}(e'))\widetilde{J}_{e'}^{\delta,\gamma}(y)$$
$$+ \gamma^3 \sum_{e',e'' \in \mathcal{V}} \mathbb{E}_\mu(\overline{\xi}(e)\overline{\xi}(e')\overline{\xi}(e''))\widetilde{J}_{e',e''}^{\delta,\gamma}(y) + O(\gamma^4),$$

where $J_e^{\delta,\gamma}(y)$ is given in Lemma 2, and where

$\widetilde{J}_{e',e''}^{\delta,\gamma}(y)$
$$= \mathbb{E}_\mu(G_{U,\delta}^{\omega^{\gamma,y}}(0,y)(\delta G_\delta^{\omega^{\gamma,y}}(y+e',y) - G_\delta^{\omega^{\gamma,y}}(y,y))$$
$$\times (\delta G_\delta^{\omega^{\gamma,y}}(y+e'',y) - G_\delta^{\omega^{\gamma,y}}(y,y)))[\mathbb{E}_\mu(G_\delta^{\omega^{\gamma,y}}(0,y))]^{-1},$$

and where, as usual, $|O(\gamma^4)| \leq C\gamma^4$, for a constant $C > 0$, depending only on $\kappa_0$, $d$. In Lemma 3, we estimated the limit of $\widetilde{J}_{e'}^{\delta,\gamma}(y)$ when $\delta$ goes to 1, by $J_{e'}^\gamma$, up to order 2 in $\gamma$. We can easily get an estimate of $\widetilde{J}_{e',e''}^{\delta,\gamma}$ at order 1 in $\gamma$, simply by expanding the terms $G^{\omega^{\gamma,y}}(y+e',y) - G^{\omega^{\gamma,y}}(y,y)$ and $G^{\omega^{\gamma,y}}(y+e'',y) - G^{\omega^{\gamma,y}}(y,y)$ at the point $G^{p_0}(y+e',y) - G^{p_0}(y,y)$ and $G^{p_0}(y+e'',y) - G^{p_0}(y,y)$, and using the transformation of step 2 to bound the rest. □

LABORATOIRE DE PROBABILITÉS
ET MODÈLES ALÉATOIRES
UNIVERSITÉ PARIS 6
4 PLACE JUSSIEU
75252 PARIS CEDEX 5
FRANCE

ECOLE NORMALE SUPÉRIEURE
DÉPARTEMENT DE MATHÉMATIQUES
ET APPLICATIONS
45 RUE D'ULM
75005 PARIS
FRANCE
E-MAIL: sabot@ccr.jussieu.fr
URL: www.proba.jussieu.fr/pageperso/sabot/index.html